\documentclass[10pt,draft]{amsart}
\usepackage{amsmath,amssymb,amsthm}
\begin{document}
\newtheorem{lem}{Lemma}[section]
\newtheorem{prop}{Proposition}[section]
\newtheorem{cor}{Corollary}[section]
\numberwithin{equation}{section}
\newtheorem{thm}{Theorem}[section]
\theoremstyle{remark}
\newtheorem{example}{Example}[section]
\newtheorem*{ack}{Acknowledgment}
\theoremstyle{definition}
\newtheorem{definition}{Definition}[section]
\theoremstyle{remark}
\newtheorem*{notation}{Notation}
\theoremstyle{remark}
\newtheorem{remark}{Remark}[section]
\newenvironment{Abstract}
{\begin{center}\textbf{\footnotesize{Abstract}}%
\end{center} \begin{quote}\begin{footnotesize}}
{\end{footnotesize}\end{quote}\bigskip}
\newenvironment{nome}
{\begin{center}\textbf{{}}%
\end{center} \begin{quote}\end{quote}\bigskip}

\newcommand{\triple}[1]{{|\!|\!|#1|\!|\!|}}
\newcommand{\xx}{\langle x\rangle}
\newcommand{\ep}{\varepsilon}
\newcommand{\al}{\alpha}
\newcommand{\be}{\beta}
\newcommand{\de}{\partial}
\newcommand{\la}{\lambda}
\newcommand{\La}{\Lambda}
\newcommand{\ga}{\gamma}
\newcommand{\del}{\delta}
\newcommand{\Del}{\Delta}
\newcommand{\sig}{\sigma}
\newcommand{\ome}{\omega}
\newcommand{\Ome}{\Omega}
\newcommand{\C}{{\mathbb C}}
\newcommand{\N}{{\mathbb N}}
\newcommand{\Z}{{\mathbb Z}}
\newcommand{\R}{{\mathbb R}}
\newcommand{\Rn}{{\mathbb R}^{n}}
\newcommand{\Rnu}{{\mathbb R}^{n+1}_{+}}
\newcommand{\Cn}{{\mathbb C}^{n}}
\newcommand{\spt}{\,\mathrm{supp}\,}
\newcommand{\Lin}{\mathcal{L}}
\newcommand{\SSS}{\mathcal{S}}
\newcommand{\F}{\mathcal{F}}
\newcommand{\xxi}{\langle\xi\rangle}
\newcommand{\eei}{\langle\eta\rangle}
\newcommand{\xei}{\langle\xi-\eta\rangle}
\newcommand{\yy}{\langle y\rangle}
\newcommand{\dint}{\int\!\!\int}
\newcommand{\hatp}{\widehat\psi}
\renewcommand{\Re}{\;\mathrm{Re}\;}
\renewcommand{\Im}{\;\mathrm{Im}\;}

\title[On the decay of solutions to a class of defocusing NLS]%
{{On the decay of solutions to a class of defocusing NLS}}
\author{}




\author{Nicola Visciglia}

\address{Nicola Visciglia\\
Dipartimento di Matematica Universit\`a di Pisa\\
Largo B. Pontecorvo 5, 56100 Pisa, Italy}

\email{viscigli@dm.unipi.it}

\maketitle
\date{}

\begin{abstract}
We consider the following family of Cauchy problems: 
\begin{equation*}
i\partial_t u= \Delta u - u|u|^\alpha, (t,x) \in \R \times \R^d
\end{equation*}
$$u(0)=\varphi\in H^1(\R^d)$$
where
$0<\alpha<\frac 4{d-2}$ for $d\geq 3$
and $0<\alpha<\infty$ for $d=1,2$.
We prove
that the $L^r$-norms of the 
solutions decay as $t\rightarrow \pm \infty$, provided that
$2<r<\frac{2d}{d-2}$ when $d\geq 3$ and $2<r<\infty$ when $d=1,2$.
In particular we extend previous results obtained 
in \cite{GinibreVelo} for $d\geq 3$
and in \cite{Nakanishi} for $d=1,2$, 
where the same decay results are proved 
under the extra assumption $\alpha >\frac 4d$.
\end{abstract}

\vspace{2cm}

\noindent This paper is devoted to the analysis of 
some asymptotic properties of solutions 
to the following family
of defocusing NLS:
\begin{equation}\label{NLS}
i\partial_t u= \Delta u - u|u|^\alpha, (t,x) \in \R \times \R^d
\end{equation}
$$u(0)=\varphi\in H^1(\R^d)$$
where
\begin{equation}\label{alpha}
0<\alpha<\frac 4{d-2} \hbox{ for } d\geq 3
\hbox{ and } 0<\alpha<\infty \hbox{ for } d=1,2.
\end{equation}
A lot of attention has been devoted in the literature 
to the Cauchy problem \eqref{NLS}. In particular 
the questions of local and global well-posedness
and scattering theory have been extensively studied.
There exists an huge literature on the field and for this reason
we cannot be exhaustive in the bibliography,
however for the moment we would like to quote
the book \cite{Cazenave}
for an extended description of the topics mentioned above
and also for an extended bibliography.

It is well--known from \cite{GinibreVelo0} 
(see also \cite{Kato} for the more general
question of unconditional uniqueness)
that, under the assumptions
\eqref{alpha} on $\alpha$, for every initial data
$\varphi\in H^1(\R^d)$ there exists
a unique global solution $u(t, x)\in {\mathcal C}
(\R; H^1(\R^d))$ to \eqref{NLS}.  
Moreover the global solutions of \eqref{NLS} satisfy
the following conservation laws:
\begin{equation}\label{charge}
\|u(t,x)\|_{L^2(\R^d)}\equiv const \hbox{ } \forall t\in\R
\end{equation}
and
\begin{equation}\label{energy}
\frac 12 \|\nabla_x u(t,x)\|^2_{L^2(\R^d)} + \frac{1}{\alpha+2}
\|u(t,x)\|_{L^{\alpha+2}(\R^d)}^{\alpha+2}\equiv const \hbox{ }
\forall t\in \R.
\end{equation} 
The main contribution of this paper concerns the decay, 
in suitable Lebesgue spaces, of 
the global solutions to \eqref{NLS}
as $t\rightarrow \pm \infty$.

\begin{thm}\label{main}
Assume $\alpha$ as in \eqref{alpha}.
Let 
$u(t,x)\in {\mathcal C}(\R;H^1(\R^d))$ be the unique global 
solution to
$$i\partial_t u= \Delta u - u|u|^\alpha, (t,x) \in \R \times \R^d$$
$$u(0)=\varphi\in H^1(\R^d).$$
Then for every $2<r<\frac{2d}{d-2}$ when $d\geq 3$
and for every $2<r<\infty$ when $d=1,2$, we have: 
\begin{equation}\label{decay}
\lim_{t\rightarrow \pm \infty }\|u(t,x)\|_{L^r(\R^d)}= 0.
\end{equation}
Moreover in the case $d=1$ we also have:
\begin{equation}\label{extra}
\lim_{t\rightarrow \pm \infty} \|u(t, x)\|_{L^\infty(\R)}
=0.
\end{equation}
\end{thm}

\begin{remark}
The proof of \eqref{extra} follows easily by combining
the conservation law \eqref{energy}, \eqref{decay}
and the following Gagliardo-Nirenberg inequality:
$$\|u(t,x)\|_{L^\infty(\R)}\leq C \|\partial_x u(t,x)\|_{L^2(\R)}^\frac 14 
\|u(t, x)\|_{L^6(\R)}^\frac 34.$$
Hence we shall focus in the sequel on the proof
of \eqref{decay}.
\end{remark}

\begin{remark}
Let us underline that the original proof of Theorem \ref{main}
in the case $\frac 4d<\alpha<\frac 4{d-2}$
and $d\geq 3$ is given in \cite{GinibreVelo}. 
In fact this is the basic step on which 
the scattering theory
in the energy space $H^1(\R^n)$ is based.
More precisely  
once the decay
of some $L^r$ norm for solutions to \eqref{NLS}
is known and $\frac 4d <\alpha <\frac4{d-2}$, then the estimates in Strichartz spaces  
follow almost
immediately, and in turn this implies easily 
the asymptotic completeness.
Hence the main novelty in our result
is that we prove dispersion of solution
to NLS also in the case $0<\alpha\leq \frac{4}{d}$.\\ 
In the case $d=1,2$ and $\frac 4n < \alpha <\infty$ the content 
of Theorem \ref{main} can be deduced from \cite{Nakanishi}.
However we point out that also in dimensions $d=1,2$ Theorem \ref{main}
covers the range
$0<\alpha\leq \frac 4d$.
\end{remark}

\begin{remark}
Theorem \ref{main} could be proved 
by the conformal conservation law 
provided that the initial data $\varphi$ belongs
to suitable weighted $L^2$ spaces. 
In fact in this case it can be deduced also a decay rate of the solution
(see Theorem 7.3.1 in \cite{Cazenave}).
However we emphasize that in general the decay of solutions
to \eqref{NLS} with initial data in $H^1(\R^d)$ and $0<\alpha\leq \frac 4n$
was a completelely open question.
\end{remark}

\begin{remark}
The proof of Theorem \ref{main}
follows by a combination of Strichartz estimates with 
the Interaction Morawetz Estimates (see \cite{CKSTT}
for the dimension $d\geq 3$ and 
\cite{CollianderGrillakis}, \cite{PlanchonVega}
for the dimensions $d=1,2$).
In particular, as a consequence of Theorem \ref{main}, 
we provide a new proof of the scattering
results in \cite{GinibreVelo} and \cite{Nakanishi}.
\end{remark}

\begin{remark}\label{prev}
Actually one can show a slightly stronger version of Theorem \ref{main}.
More precisely in the case $d\geq 3$ and $\frac 4d <\alpha
<\frac 4{d-2}$ we also have
\begin{equation}\label{impro}
\lim_{t\rightarrow \pm \infty} \|u(t, x)\|_{L^{\frac{2d}{d-2}}(\R^d)}=0.
\end{equation}
The proof of this fact goes as follows. 
Once \eqref{decay} is proved and $\alpha$ is as above, then
by using classical arguments (see \cite{Cazenave}), we can
construct the scattering operator in the energy space. More precisely
given any $\varphi \in H^1(\R^d)$
there exist $\varphi_\pm\in H^1(\R^d)$ such that
\begin{equation}\label{free}
\lim_{t\rightarrow \pm \infty}
\|u(t,x) - e^{it\Delta} \varphi_\pm\|_{H^1(\R^n)}=0,
\end{equation}
where $u(t, x)$ denotes the corresponding solution to \eqref{NLS}.
Hence we get
\begin{equation}\label{mink}
\|u(t, x)\|_{L^{\frac{2d}{d-2}}(\R^d)}
\leq \|u(t,x) - e^{it\Delta} \varphi_\pm \|_{L^{\frac{2d}{d-2}}(\R^d)}
+ \|e^{it\Delta} \varphi_\pm \|_{L^{\frac{2d}{d-2}}(\R^d)}\end{equation}
$$ \leq C \|u(t,x) - e^{it\Delta} \varphi_\pm \|_{H^1(\R^d)}
+ \|e^{it\Delta} \varphi_\pm \|_{L^{\frac{2d}{d-2}}(\R^d)},$$
where we have used the Sobolev embedding at the last step.
On the other by combining the dispersive estimate
$$\|e^{it\Delta} \psi\|_{L^{\frac{2d}{d-2}}(\R^d)}\leq \frac C{t} \|\psi\|_{L^{\frac{2d}{d+2}}(\R^d)},$$
with the Sobolev embedding and with a density argument, we can deduce
easily that
\begin{equation}\label{freeop}
\lim_{t\rightarrow \pm \infty}
\|e^{it\Delta} \psi\|_{L^{\frac{2d}{d-2}}(\R^d)}=0 \hbox{ } \forall \psi \in H^1(\R^d).
\end{equation}
Hence by combining \eqref{free}, \eqref{mink} and \eqref{freeop}
we get \eqref{impro}.
\end{remark}

\begin{remark}
Arguing as in remark \ref{prev}
and by using the embedding
$L^\infty(\R^2)\subset BMO(\R^2)$ 
and $H^1(\R^2)\subset BMO(\R^2)$, we can deduce:
$$\lim_{t \rightarrow \pm \infty}
\|u(t,x)\|_{BMO(\R^2)}=0,$$
provided that $u(t, x)$ solves \eqref{NLS} with $d=2$ and $\alpha>2$.
\end{remark}

\noindent Along this paper a fundamental role will be played by 
the Strichartz estimates for the propagator
$e^{it\Delta}$, hence for the sake of completeness they will be stated below.
\\
First we need to introduce some notations that will be useful in the sequel.
\\ 
\\\noindent For any subinterval $I\equiv (0, T)$ of $\R$  
and for every $p,q\in [0, \infty]$ we define the mixed space-time norms
\begin{equation}\label{eq.spt}
    \|u\|_{L^p_TL^q_x}\equiv \left(\int_{I}\|u(t,\cdot)\|_{L^{q}(\R^{n})}^{p}dt \right)^{1/p}.
    \end{equation}
In the case $I\equiv \R$ we write
$\|u\|_{L^p_t L^q_x}$.     \\
\noindent We shall also use the notation
$\|\varphi\|_{L^r_x}\equiv \|\varphi\|_{L^r(\R^d)}$ for every $1\leq r\leq \infty$
and $\|\varphi\|_{H^1_x}\equiv \|\varphi\|_{H^1(\R^d)}$.     
\\\noindent Given $d\geq 1$ we say that the pair $(p,q)$ is $d$-\emph{(Schr\"odinger) admissible} if
\begin{equation}\label{q.strich}
    \frac 2 p + \frac d {q} = \frac d 2,\qquad p,q\in[2,\infty],\qquad
    (d,p,q)\neq(2,2,\infty).
\end{equation}
\\
\\
\noindent We can now state the Strichartz estimates for the free
propagator.
For any $d$-\emph{(Schr\"odinger) admissible} couples $(p,q)$ and
$(\tilde p,\tilde q)$ there exists a constant $C(p,\tilde p)$ such that, for all 
$T>0$, for all functions $u_{0}\in L^2_x$,
and $F(t,x)\in L^{\tilde p'}_TL^{\tilde q'}_x$ the following inequalities hold:
\begin{equation}\label{eq.strich}
    \biggl\| e^{it\Delta}u_{0}\;\biggr\|_{L^p_TL^q_x}\leq C(p,\tilde p)\;\|u_{0}\|_{L^{2}_x}
\end{equation}
\begin{equation}\label{eq.strichF}
    \biggl\|\int_{0}^{t}e^{i(t-s)\Delta}F(s)\;ds\biggr\|_{L^p_T L^q_x}\leq 
       C(p,\tilde p)  \;  \bigl\| F\bigr\|_{L^{\tilde p'}_TL^{\tilde q'}_x},
\end{equation}
where we have used a prime to denote conjugate indices.
Note that the constant in the previous estimates
are independent of the interval $T>0$. 
For a proof of the Strichartz estimates in the non end--point case
see
\cite{Cazenave}, for the general case see 
\cite{KT}.
\\
\\

\noindent {\bf Acknowledgement:} {\em 
the author is greateful to T. Cazenave
for bringing to his attention the question studied in this paper
and for many useful advices, and to L. Vega for interesting discussions about the interaction
Morawetz estimates.}

\section{Proof of Theorem \ref{main}}
We shall need the following
\begin{lem}\label{prel}
Let $\chi\in C^\infty_0(\R^d)$
be a cut--off function and
$\varphi_n\in H^1_x$ be a sequence such that
$M\equiv \sup_{n\in\N} 
\|\varphi_n\|_{H^1_x}<\infty$
and $\varphi_n\rightharpoonup \bar \varphi \hbox{ in } H^1_x$.
Let $u_n(t, x), \bar u(t,x)\in {\mathcal C}(\R; H^1_x)$ 
be the corresponding solutions
to \eqref{NLS} with initial data $\varphi_n$
and $\bar \varphi$ respectively.
Then
for every $\epsilon>0$ there exists
$T(\epsilon)>0$ and $\nu (\epsilon) \in \N$ such that 
$$\sup_{t\in (0, T(\epsilon))}\|\chi(x) (u_n(t,x) - \bar u(t,x))\|_{L^2_x}\leq \epsilon
\hbox{ } \forall n>\nu(\epsilon).$$
\end{lem}

\noindent{\bf Proof.}
By combining \eqref{charge} and \eqref{energy}
it is easy to deduce that
\begin{equation}\label{bound}
\sup_{n\in \N,t\in \R} 
\{\|u_n(t, x)\|_{H^1_x}, \| \bar u(t, x)\|_{H^1_x}\}<\infty.
\end{equation}
By using the Rellich compactness theorem
we have
\begin{equation}\label{rellich}
\lim_{n\rightarrow \infty}
\|\chi(x)(\varphi_n - \bar \varphi)\|_{L^2_x}=0.
\end{equation}
Next we introduce the functions 
$$v_n(t, x)\equiv \chi(x) u_n(t, x) 
\hbox{ and } \bar v(t,x)\equiv \chi(x) \bar u(t, x)$$
that solve the following
Cauchy problems:
$$i \partial_t v_n= \Delta v_n
-2 \nabla \chi \cdot \nabla u_n - u_n\Delta \chi  - \chi u_n|u_n|^\alpha$$
$$v_n(0)= \chi(x) \varphi_n$$
and 
$$i \partial_t \bar v= \Delta \bar v
-2 \nabla \chi \nabla \cdot \bar u - \bar u\Delta \chi  - \chi \bar u|\bar u|^\alpha$$
$$\bar v(0)= \chi(x) \bar \varphi.$$
By using the integral formulation of the previous Cauchy problems
we deduce:
\begin{equation}\label{integral}
v_n(t,x)-\bar v(t,x)= e^{it\Delta} [\chi(x) (\varphi_n-\bar \varphi)]
\end{equation}
$$+ i \int_0^t e^{i(t-s)\Delta}[2 \nabla \chi \cdot \nabla (u_n(s)-\bar u(s)) +
 (u_n(s)-\bar u(s)) \Delta \chi] dx  
$$$$ + i \int_0^t e^{i(t-s)\Delta}
 [\chi(x) (u_n(s)|u_n(s)|^\alpha -  \bar u(s)|\bar u(s)|^\alpha)] ds.
$$

\noindent Next we split the proof in two cases.
\\

{\em First case: $d\geq 3$}
\\

\noindent We fix the following $d$-\emph{(Schr\"odinger) admissible} couple 
$$(p(d), q(d))\equiv \left  (\frac 8{\alpha(d-2)},\frac{4d}{2d-\alpha d + 2\alpha}
\right )$$ and by using the estimates 
\eqref{eq.strich} and \eqref{eq.strichF}
(where we make the choices $(\tilde p, \tilde q)=(p(d), q(d))$
and $(\tilde p, \tilde q)=(\infty,2)$) we get:
$$\|v_n - \bar v\|_{L^{p(d)}_T L^{q(d)}_x}
\leq C \| \chi(x) (\varphi_n - \bar \varphi)\|_{L^2_x}
+ C \|\nabla \chi \cdot \nabla (u_n- \bar u) \|_{L^1_TL^2_x}
$$$$+ C \| (u_n- \bar u) \Delta \chi\|_{L^1_T L^2_x} 
+ C \| \chi (u_n|u_n|^\alpha - \bar u |\bar u|)^\alpha\|_{L^{p(d)'}_TL^{q(d)'}_x},$$
that in conjunction with \eqref{bound} and  
with the H\"older inequality implies:
$$\|v_n - \bar v\|_{L^{p(d)}_T L^{q(d)}_x}
\leq C \| \chi(x) (\varphi_n - \bar \varphi) \|_{L^2_x}
+ C T 
$$$$+ C \| \chi(x) (u_n - \bar u)\|_{L^{p(d)'}_T L^{q(d)}_x}
\sup_{t\in (0, T)} 
\left (\|u_n(t)\|^\alpha_{L^{\frac{2d}{d-2}}_x} + \|\bar u(t)\|^\alpha_{L^{\frac{2d}{d-2}}_x}\right ).$$
By using now the Sobolev embedding 
$H^1_x\subset L^\frac{2d}{d-2}_x$, \eqref{bound} and the H\"older inequality
in the time variable, we get:
$$\|v_n - \bar v\|_{L^{p(d)}_T L^{q(d)}_x}
\leq C \| \chi(x) (\varphi_n - \bar \varphi)\|_{L^2_x}
+ C T 
$$$$+ C  T^\frac{p(d)-2}{p(d)} \| v_n - \bar v\|_{L^{p(d)}_T L^{q(d)}_x}.$$
By combining this estimate with \eqref{rellich} we deduce that
for every $\epsilon>0$ there exist $\nu(\epsilon)$ and $T(\epsilon)$
such that
\begin{equation}\label{almost}
\|v_n - \bar v\|_{L^{p(d)}_{T(\epsilon)} L^{q(d)}_x}\leq \epsilon
\hbox{ } \forall n>\nu(\epsilon).\end{equation}
Next we consider again \eqref{integral}
and we use again the Strichartz estimates 
with the choice $(p, q)=(\infty, 2)$ and $(\tilde p, \tilde q)$ as above,
and arguing as above
we deduce:
$$\|v_n - \bar v\|_{L^\infty_{T(\epsilon)}L^2_x}
\leq C \| \chi(x) (\varphi_n - \bar \varphi)\|_{L^2_x}
+ C T(\epsilon)
$$$$+ C T(\epsilon)^\frac{p(d)-2}{p(d)} \| v_n - \bar v\|_{L^{p(d)}_{T(\epsilon)} 
L^{q(d)}_x}.$$
By combining this estimate with \eqref{rellich}
and \eqref{almost}
we deduce the desired result.
\\

{\em Second case: $d=1,2$}
\\

\noindent The proof is similar to the case $d\geq 3$ provided that
we make respectively the following choice of
 $1$-\emph{(Schr\"odinger) admissible} and  $2$-\emph{(Schr\"odinger) admissible} 
 couples:
$$(p(1), q(1))=(\infty,2) \hbox{ and } (p(2), q(2))= (4,4).$$

\hfill$\Box$

\vspace{0.1cm}

\noindent{\bf Proof of Thm \ref{main}}
We shall prove \eqref{decay}
for $t\rightarrow \infty$ (the case $t\rightarrow -\infty$
can be treated in a similar way). 
We split the proof in two cases.
\\

{\em First case: $d\geq 3$}
\\

\noindent Notice that by combining \eqref{charge} and \eqref{energy}
with the H\"older inequality, it is enough to prove \eqref{decay}
for $r=\frac{2d+4}{d}$.
Next we recall the following consequence
of the Gagliardo--Nirenberg inequality for $d\geq 3$:
\begin{equation}\label{GN}
\|\psi\|_{L^\frac{2d+4}{d}_x}^\frac{2d+4}d\le 
C \left (\sup_{x\in \R^d}\|\psi\|_{L^2(Q_1(x))}\right )^{\frac 4d}  
 \| \psi\|_{H^1_x}^{2},
\end{equation}
where $Q_r(x)$ denote
the cube in $\R^d$ centered in $x$ whose edge
has lenght $r$.
Moreover as a consequence of
\eqref{charge} and \eqref{energy}
we get 
\begin{equation}\label{apb}
\sup_{t\in\R} \|u(t, x)\|_{H^1_x}<\infty.
\end{equation}
\noindent Next we assume by the absurd that there is a sequence $t_n\rightarrow \infty$ such that
$$\|u(t_n,x)\|_{L^\frac{2d+4}d_x}\geq \epsilon_0>0.$$
Then by combining \eqref{apb} with \eqref{GN}, where we choose $\psi\equiv u(t_n, x)$,
we deduce the existence of a sequence
$x_n \in \R^d$ such that
\begin{equation}\label{positive}
\|u(t_n,x)\|_{L^2(Q_1(x_n))}\geq \epsilon_1>0.
\end{equation}
Next we introduce the functions
\begin{equation}\label{psi}
\varphi_n(x)\equiv u(t_n, x+x_n).\end{equation}
which are bounded in $H^1_x$ by \eqref{apb}.
By combining the compactness of the Sobolev embedding on bounded set
with 
\eqref{positive}, we have that
up to subsequence $\varphi_n$ converges
weakly in $H^1_x$ to a nontrivial function $\bar \varphi
\in H^1_x$ such that
\begin{equation}\label{bubble}
\|\bar \varphi\|_{L^2(Q_1(0))}\geq \epsilon_1>0.
\end{equation}
We are now in condition to apply Lemma \ref{prel}
where we choose the function $\chi(x)$ as any cut--off function
supported in $Q_2(0)$ and
such that $\chi(x)\equiv 1$ on the cube $Q_1(0)$.
We also introduce the functions $u_n(t, x)$ and $\bar u(t,x)$
as the solutions to:
$$i\partial_t u_n= \Delta u_n - u_n|u_n|^\alpha, (t,x) \in \R \times \R^d$$
$$u_n(0)=\varphi_n\in H^1_x$$
and
$$i\partial_t \bar u= \Delta \bar u - \bar u|\bar u|^\alpha, (t,x) \in \R \times \R^d$$
$$\bar u(0)=\bar \varphi\in H^1_x.$$
Notice that by combining \eqref{bubble} 
with a continuity argument we deduce the existence of
$\bar T>0$ such that
\begin{equation}\label{almo}
\inf_{t\in (0, \bar T)}
\|\chi(x) \bar u(t,x)\|_{L^2_x}\geq \frac{\epsilon_1}2>0.
\end{equation}
By Lemma \ref{prel} there exist
$\tilde T>0$ and $\nu\in \N$ such that
\begin{equation}\label{alm}
\sup_{t\in (0, \tilde T)}\|\chi(x) (u_n(t) - \bar u(t))\|_{L^2_x}\leq \frac{\epsilon_1} 4
\hbox{ } \forall n>\nu,
\end{equation}
and by combining \eqref{almo} and \eqref{alm}
we get
\begin{equation}\label{abs}
\|\chi(x) u_n(t,x)\|_{L^2_x}
\geq - \|\chi(x) (u_n(t,x) - \bar u(t,x))\|_{L^2_x}
+ \|\chi(x) \bar u(t,x)\|_{L^2_x}\geq \frac {\epsilon_1}4
\end{equation}$$
\hbox{ } \forall t<T_0\equiv \min \{\bar T, \tilde T\}, \forall n>\nu.$$
Notice that, due to the properties of $\chi(x)$, \eqref{abs} implies
\begin{equation}
\|u_n(t, x)\|_{L^2(Q_2(0))}\geq \frac {\epsilon_1} 4 
\hbox{ } \forall t<T_0, \forall n>\nu.
\end{equation}
On the other hand by the translation invariance of
NLS and due to the definition on $\varphi_n$ (see \eqref{psi}),
it is easy to deduce that the previous 
estimate is equivalent to the following one:
\begin{equation}\label{L2}
\|u(t, x)\|_{L^2(Q_2(x_n))}\geq \frac {\epsilon_1} 4
\hbox{ } \forall t\in (t_n,t_n+T_0), \forall n>\nu.
\end{equation}
Next we show that \eqref{L2} lead to a contradiction
and it will complete the proof of \eqref{decay}
when $d\geq 3$.\\

\noindent Choose $a(x)\equiv |x|$
in the inequality written at page
10 in \cite{CKSTT}.
Then this implies 
\begin{equation}\label{CKS}
\int_\R \int\int_{\R^d_x\times \R^d_y}\frac{|u(t,x)|^2 |u(t,y)|^2}{|x-y|^3} dxdydt<\infty
\hbox{ for } d\geq 4
\end{equation}
and
\begin{equation}\label{CKS1}
\|u(t,x)\|_{L^4_{t,x}}<\infty \hbox{ for } d=3.
\end{equation}
Notice that since $t_n\rightarrow \infty$ we can assume
up to subsequence 
that the sets $(t_n, t_n + T_0)$ are disjoint sets.
Next we use \eqref{L2} and we get
$$\int_\R \int\int_{\R^d_x\times \R^d_y}
\frac{|u(t, x)|^2 |u(t, y)|^{2}}{|x-y|^3} dtdxdy$$$$
\geq C \sum_{n\in N} \int_{t_n}^{t_n+T_0} \int \int_{Q_2(x_n)\times Q_2(x_n)}
|u(t, x)|^2 |u(t, y)|^{2}| dtdxdy =\infty$$
and this is in contradiction with \eqref{CKS}.
This complete the proof of \eqref{decay}
for $d\geq 4$. Similarly for $d=3$ we deduce by \eqref{L2}
that $$\|u(t, x)\|_{L^4(Q_2(x_n)
\times (t_n, t_n+ T_0))}^4
\geq C \epsilon_1^4 T_0 $$
and also in this case we get easily a contradiction with \eqref{CKS1}.
\\

{\em Second case: $d=1,2$}
\\

\noindent As in the previous case
it is sufficient to prove 
\eqref{decay} for $r=3$. In order to do that
we shall need the following
version of \eqref{GN} in dimensions $d=1,2$:
$$\|\psi\|_{L^3_x}^3\leq C \left (\sup_{x\in \R^d}\|\psi\|_{L^2(Q_1(x))}\right ) \|\psi\|_{H^1_x}^2.$$
Arguing as in the previous case
we can deduce that if \eqref{decay} is false with $r=3$,
then there exist $(t_n,x_n)\in \R\times \R^d$ such that
$t_n\rightarrow \infty$ as $n\rightarrow \infty$
and 
\begin{equation}\label{L2d12}
\|u(t, x)\|_{L^2(Q_2(x_n))}\geq \epsilon_2>0
\hbox{ } \forall t\in (t_n,t_n+T_0).
\end{equation}
It is now easy to deduce (arguing as we did above in the case
$d=3$) that this lead to a contradiction
with the following a--priori bounds
proved in 
\cite{CollianderGrillakis}
and \cite{PlanchonVega}:
$$
\|u\|_{L^4_t L^8_x} <\infty \hbox{ for } d=2
\hbox{ and }
\|u\|_{L^{\alpha+4}_{t,x}} <\infty \hbox{ for } d=1,
$$
(more precisely the estimates above follow for instance from Theorems
1.1,1.2 in \cite{CollianderGrillakis}).
Hence the proof of \eqref{decay} is complete also for $d=1,2$.

\hfill$\Box$

\end{document}